\documentclass{article}
\usepackage{amssymb}
\usepackage{a4}
\usepackage{psfig}
\usepackage{epsfig}
\usepackage{latexsym}
\usepackage{pictex}
\usepackage{graphicx}
\usepackage{amsmath}
\usepackage{amsfonts}
\begin{document}
%%%%%%%%%%%%%%%%%%%%%%%%%%%%%%%%%%%%%%%%%%%%%%%%%%%%%%%%%%%%%%%%%%%%%%
%	spaces for your own definitions follows
%%%%%%%%%%%%%%%%%%%%%%%%%%%%%%%%%%%%%%%%%%%%%%%%%%%%%%%%%%%%%%%%%%%%%%
\renewcommand {\a}{ \alpha }
\renewcommand {\b}{\beta}
\newcommand{\D}{\Delta}
\newcommand{\e}{\epsilon}
\newcommand{\vare}{\varepsilon}
\newcommand{\g}{\gamma}
\newcommand{\G}{\Gamma}
\newcommand{\varf}{\varphi}
\renewcommand{\d}{\delta}
\newcommand{\s}{\sigma}
\renewcommand{\l}{\lambda}
\renewcommand{\t}{\theta}
\newcommand{\vart}{\vartheta}

\newcommand{\R}{ \mathbb R}
\newcommand{\N}{ \mathbb N}

\newcommand {\ba}{\mathbf a}
\newcommand {\BA}{\mathbf A}

\newcommand{\CA}{\mathcal A}
\newcommand{\CE}{\mathcal E}
\newcommand{\CF}{\mathcal F}
\newcommand{\CH}{\mathcal H}
\newcommand{\CV}{\mathcal V}

\newcommand{\plainH}[1]{\textup{{\textsf{H}}}^{#1}}
\newcommand{\plainL}[1]{\textup{{\textsf{L}}}^{#1}}

\newcommand{\CI}{\mathcal I} 
\newcommand{\CM}{\mathcal M}
\newcommand{\f}{\varphi}
\newcommand{\plainM}{\textup{\textsf{M}}}
\newcommand \gen{\text{Gen}}
\newcommand \Dom{\text{Dom}}
\newcommand\Quad{\text{Quad}}
\newcommand\prob{\text{prob}}
\newcommand\res{\restriction}

\newtheorem{defi}{Definition}[section]
\newtheorem{theo}[defi]{Theorem}%[Theorem]
\newtheorem{prop}[defi]{Proposition}
\newtheorem{lemma}[defi]{Lemma}
\newtheorem{cor}[defi]{Corollary}%[Corollary]
%

%\theoremstyle{definition}
%\newtheorem{definition}{Definition}
%\newtheorem{example}{Example}
%\newtheorem{remark}{Remark}

%%%%%%%%%%%%%%%%%%%%%%%%%%%%%%%%%%%%%%%%%%%%%%%%%%%%%%%%%%%%%%%%%%%%%%
%	Please fill in your details below,
%	address at the end of the file
%%%%%%%%%%%%%%%%%%%%%%%%%%%%%%%%%%%%%%%%%%%%%%%%%%%%%%%%%%%%%%%%%%%%%%

\author{Michael Solomyak 	}

\title{Laplace and Schr\"odinger operators on regular metric trees:
the discrete spectrum case}

\date{\em Dedicated to Prof. Hans Triebel on the occasion of his 65th
birthday}

%\date{}

\maketitle

%%%%%%%%%%%%%%%%%%%%%%%%%%%%%%%%%%%%%%%%%%%%%%%%%%%%%%%%%%%%%%%%%%%%%%%
%	You can insert an abstract here if you want. 
%%%%%%%%%%%%%%%%%%%%%%%%%%%%%%%%%%%%%%%%%%%%%%%%%%%%%%%%%%%%%%%%%%%%%%%
%\begin{abstract}

%\end{abstract}
%%%%%%%%%%%%%%%%%%%%%%%%%%%%%%%%%%%%%%%%%%%%%%%%%%%%%%%%%%%%%%%%%%%%%%%
%	Please insert the article body now
%%%%%%%%%%%%%%%%%%%%%%%%%%%%%%%%%%%%%%%%%%%%%%%%%%%%%%%%%%%%%%%%%%%%%%%

\section{Introduction} Spectral theory of differential
operators on metric trees is an interesting branch of such theory 
on general metric graphs. Among the trees, the so-called regular
trees are of particular interest due to their very special geometry. 

Let $\G$ be a tree rooted at some vertex $o$ and having infinitely many
edges.  Below
$|x|$ stands for the distance between a point $x\in\G$ and the root $o$.
For a vertex $x$,  its {\sl generation} is the number of vertices lying 
between $o$ and $x$ (including $x$ but excluding $o$). We say
that a tree $\G$ is {\sl regular} if for any vertex $x$ the quantity $|x|$ and
the number of edges
emanating from $x$ depend only on the generation of $x$; see 
Definition 2.1 for the
more
detailed description. 

The regular trees are highly symmetric. This allows one to construct an
orthogonal decomposition of the space $\plainL2(\G)$ which reduces the
Schr\"odinger operator $\BA_V=-\boldsymbol{\D}+V$ with any symmetric 
(i.e. depending on
$|x|$) potential $V$. We call this decomposition the {\sl basic
decomposition} of $\plainL2(\G)$. The parts of $\BA_V$  in the components 
of the basic decomposition are denoted by $\CA_{V,k}$. Each operator
$\CA_{V,k}$ appears with a multiplicity which rapidly grows as $k\to\infty$.
The study of the spectrum
$\s(\BA_V)$ of the operator $\BA_V$ is thus reduced to the study of the spectra
of the parts $\CA_{V,k}$. Each part
can be identified with a differential operator acting in a
weighted space $\plainL2\bigl((t_k,R),g_k\bigr)$, $k=0,1,\ldots$ where the
intervals $ (t_k,R)$ and the weight functions $g_k$ are
determined by the geometry of $\G$. In particular, the quantity
$R=R(\G)=\sup\{|x|:x\in\G\}\le\infty$ is the {\sl radius} of the tree. 

According to general spectral theory, \begin{equation*}
\s(\BA_V)=\overline{\cup_{k=0}^\infty\s(\CA_{V,k})}. \end{equation*}
However, the quantitative characteristics of $\s(\BA_V)$ can
not be obtained automatically from the corresponding characteristics for
the operators $\CA_{V,k}$, due to the growing multiplicities of 
$\CA_{V,k}$ as parts of $\BA_V$.

\bigskip

The main purpose of this paper is to study the situations when the
spectrum of the Schr\"odinger operator $\BA_V$ is discrete. More exactly,
we consider two typical cases: for the regular trees of finite radius we
show that the classical Weyl asymptotic law holds for the eigenvalues of
each operator $\CA_{V,k}$ with any bounded potential $V$,
including the basic case $V=0$. If the tree $\G$ 
has finite total length, we
show that the Weyl formula holds for the whole operator $\BA_V$.  

Another case is the operator $\BA_V$ with growing potential
on a tree of infinite radius.  Under certain
assumptions about the geometry of the tree and the behaviour of the
potential, we find a version of the Weyl formula for such operators. 

To make the general picture more complete, we also present, without proofs,
some known results concerning the structure of $\s(\BA_V)$ in the case
when the spectrum is not discrete.

\bigskip

There are several papers, devoted to differential operators on
regular metric trees. In \cite{C1} the case of the {\sl homogeneous trees
$\G_b$} was considered. A regular tree is called {\sl homogeneous} if all
its edges have equal length (say $1$) and all vertices have the same
number of edges (say $b$) emanating from them. In \cite{C1} the 
potential $V$ was supposed to be periodic and even. It was shown that the
spectrum $\s(\BA_V)$ has the band-gap structure, with no more
than one eigenvalue in each gap. Actually, this eigenvalue is of infinite
multiplicity, but the multiplicities were not discussed in \cite{C1}. 
\vskip0.2cm
In \cite{NS1} the weighted spectral problems of the form
$-\boldsymbol{\D}f=\l Vf$ on
general (not necessarily regular) trees were investigated.  The estimates
and, under some additional assumptions, the asymptotic behaviour of the
eigenvalues were found. For the regular trees, 
the basic decomposition of $\plainL2(\G)$ 
was discovered in \cite{NS1}, and much more advanced results for such 
trees were obtained with the
help of this decomposition. 

\vskip0.2cm In the paper \cite{C2} the basic decomposition was
re-discovered, and a detailed spectral analysis of the operators
$\CA_{V,k}$ for the regular trees $\G$ with finite radius
was given. In particular, it was shown that if $R(\G)<\infty$
but the total length of $\G$ is infinite, the operators $\CA_{V,k}$ 
do not require the boundary condition at $t=R$. It was also shown that
each operator $\CA_{V,k}$ has compact resolvent, and the corresponding
eigenvalue distribution function $N(\l;\CA_{V,k})$ grows not faster than 
$O(\l^{1/2+\e})$ for any $\e>0$. Our Theorem 5.3 considerably refines 
this result.
 
\vskip0.2cm The main topic of \cite{NS2} is the
Hardy-type inequalities on regular trees. As a consequence, a necessary
and sufficient condition of the positive definiteness of the Laplacian was
established, see Theorem 5.6.

\vskip0.2cm In \cite{SS} operators $\BA_V$
with decaying symmetric potentials on the homogeneous trees $\G_b$ were
investigated. For $V=0$ (that is, for the free Laplacian) the spectrum was
explicitly calculated. This result is presented here as Theorem 5.7.
In accordance with the results of \cite{C1}, the spectrum has infinitely
many gaps. Perturbation of the operator $\BA_0$ by a decaying
potential may create eigenvalues in each gap, and in \cite{SS} their
behaviour was investigated in detail both for positive and
negative perturbations. 
\vskip0.2cm
It is necessary to mention here also the papers \cite{EH} and \cite{EHL},
though formally they do not deal with the Laplacian on trees. In \cite{EH}
the Hardy-type integral operators on trees were introduced in connection
with the spectral analysis of the Neumann Laplacian in 
certain irregular domains. For these ``ridged'' domains, there exists a tree
that serves as a ``ridge'', or a ``skeleton'', for the domain.
There is a close relation between the approximation numbers of the 
Hardy-type integral operators and the eigenvalues of the problem 
$-\boldsymbol{\D}f=\l Vf$ on the tree. It was the paper \cite{EH}
which attracted the author's attention to operators on trees.

In \cite{EHL} the behaviour of the approximation numbers of the
Hardy-type integral operators on trees was studied in detail, not only
in $\plainL2$-case but also in the general $\textsf L^p$-case. When applied to
the Laplacian, the estimates obtained substantially refine some results
of \cite{NS1}, Sect. 4. The asymptotic formulae found in \cite{EHL} may
provide an alternative approach to the proof of our Theorem \ref{sp:4}, (ii).

\bigskip

The structure of the present paper is as follows. The short sections 2 -- 4
contain the necessary preliminary material: the definitions of 
regular and homogeneous trees and of the Laplace and  Schr\"odinger
operators on them, and the description of the basic orthogonal decomposition
of $\plainL2(\G)$. Sect. 5 contains formulations of the main results,
the proofs are given in Sect. 6. In the final Sect. 7 we discuss the extension
of the presented results to the case of regular trees without boundary.

\section{Regular rooted trees } 
\subsection{Geometry of a tree} Let $\G$ be a rooted tree with  
the root $o$, the set of vertices $\CV=\CV(\G)$ 
and the set of edges $\CE=\CE(\G)$. We suppose that $\#\CV=\#\CE=\infty$.
Unlike the combinatorial trees, whose edges are just pairs of vertices,
each edge $e$ of a metric tree is viewed as a non-degenerate line segment. 
The distance $\rho(x,y)$ between
any two points $x,y\in\G$ (and thus the metric topology on $\G$) is 
introduced in a natural way. As it was already said in Introduction, 
$|x|$ stands
for $\rho(x,o)$. 

We write $y\preceq z$ if $|z|=|y|+\rho(y,z)$, and 
$y\prec z$ means that $y\preceq z$ and $z\neq y$.  The 
relation $\prec$ defines on $\G$ a partial ordering. If $y\prec z$, 
we denote 
\begin{equation*}
\langle y,z\rangle := \{ x\in\G:\;y\preceq x\preceq z \}.
\end{equation*}
%In particular, if $e=\langle y,z\rangle$ is an edge, we 
%say that $e$ emanates from $y$ and terminates at $z$. Only one edge, say
%$e_z^-$,
%terminates at each vertex $z\neq o$, and no edge terminates at $o$.

For any vertex $y$ its {\it{generation}} $\gen(y)$ is defined as 
\begin{equation*}
\gen(y)=\#\{x\in\CV:o\prec x\preceq y\}.
\end{equation*}
We assume that $\gen(y)<\infty$ for any vertex $y$. 
For an edge $e$ we define $\gen(e)$ as
the generation of its initial point. 

The {\sl{branching number}} $b(y)$ of a vertex $y$ is defined as 
the number of edges
emanating from $y$. We assume that $b(o)=1$ and $b(y)>1$ for $y\neq o$.
We denote by $e_y^-$ the only edge which terminates at a vertex $y\neq o$,
and by $e_y^1,\ldots,e_y^{b(y)}$ the edges emanating from any vertex
$y\in\CV$.

\begin{defi}\label{tr:def} 
We call a tree $\G$ {\sl{regular}} if all the vertices of the same generation 
have equal
branching numbers, and all the edges of the same generation are of the same 
length.
\end{defi}

In this paper we consider only the regular trees. Evidently, any such 
tree is fully determined by specifying two
number sequences, $\{b_k\}$ and $\{t_k\}$, $k=0,1,\ldots$ such that 
\begin{equation*}
b(y)=b_{\gen(y)},\ |y|=t_{\gen(y)}\qquad 
{\text{for each}}\;y\in\CV(\G).
\end{equation*}
According to our assumptions, $b_0=1$ and $b_k\ge2$ for any $k>0$.
It is 
clear that $t_0=0$ and the sequence $\{t_k\}$ is strictly increasing,
and we denote
\begin{equation*}%\label{1:rad}
R=R(\G)=\lim_{k\to\infty}t_k =\sup_{x\in\G}|x|.
\end{equation*}
We call $R(\G)$ the {\sl radius } of the tree. Another important
characteristic of a tree is its {\sl total length} (in other terminology,
{\sl volume}) 
\begin{equation*}%\label{1:tot}
|\G|=\sum_{e\in\CE(\G)}|e|.
\end{equation*}
The natural measure $dx$ on 
$\G$ is induced by the Lebesgue measure on the edges. Below
$\plainL2(\G)=\plainL2(\G,dx)$.

\subsection{Homogeneous trees} A rooted tree is called homogeneous if its 
edges are all of the 
same length
(for definiteness, of the length one) and all the vertices 
$y\neq o$ have the same branching number $b$.
A homogeneous tree is evidently regular. It is fully
determined by specifying the parameter $b$ and we use for it
the notation $\G_b$. For the tree $\G_b$ one has $t_k=k,\ k=0,1,\dots$,
and $b_k=b,\ k=1,2,\dots$. 

\section{The Laplace and the Schr\"odinger operators on a regular tree } 
The notion of differential operator on any metric graph, in particular
on a tree, is well known.
Still, for the sake of completeness we present here the variational definitions
of the Laplacian and of the Schr\"odinger operator on a tree.

We say that a scalar-valued function $f$ on $\G$ belongs to the 
Sobolev space $\plainH1=\plainH1(\G)$ if $f$ is continuous,
$f\res e\in \plainH1(e)$ for each edge $e$, and
\begin{equation}\label{h}
\| f\|^2_{\plainH1}:=\int_\G\bigl(|f'(x)|^2+|f(x)|^2\bigr)dx<\infty.
\end{equation}
The derivative of a function $f\res e$ at an interior point $x\in e$ is
always taken in the direction compatible with the 
partial ordering on $\G$. This
agreement is indifferent for the definition \eqref{h} but we shall use it
later. 

The set
$\plainH1_{\text b}$ of all boundedly supported functions 
$u\in \plainH1$
is dense in $\plainH1$. Indeed,
for any number $L>0$ let $\f_L(t)$ be the continuous 
function on $\R_+$,
which is $1$ for $t\le L$, is $0$ for $t\ge L+1$ and is linear on 
$[L,L+1]$.
Given a function $f\in \plainH1(\G)$, denote $f_L(x)=\f_L(|x|)f(x)$. 
Then
$f_L\in \plainH1_{\text b}$ and an elementary calculation shows that 
$f_L\to f$ in $\plainH1(\G)$ as $L\to\infty$. 

Along with $\plainH1(\G)$, let us introduce also its subspace of 
codimension one:
\begin{equation*}
\plainH{1,0}:=\plainH{1,0}(\G)=\bigl\{f\in \plainH{1}(\G):f(o)=0\bigr\}.
\end{equation*}
We define the Dirichlet
Laplacian $-\boldsymbol{\D}$ on $\G$ as the self-adjoint
operator in
$\plainL2(\G)$, associated with the
quadratic form $\int_\G|f'|^2dx$ considered
on the form domain $\Quad(-\boldsymbol{\D})=\plainH{1,0}(\G)$. 
It is easy to
describe the operator domain
$\Dom(\boldsymbol{\D})$ and the action of $\boldsymbol{\D}$.
Evidently $f\in\Dom(\boldsymbol{\D})\Rightarrow f\res e\in \plainH2(e)$
for each edge $e$ and the Euler -- Lagrange
equation reduces on $e$ to
$\boldsymbol{\D} f=f''$. At the root we have the boundary condition $f(o)=0$,
since $\Dom(\boldsymbol{\D})\subset \plainH{1,0}(\G)$. 
At each vertex $y\neq o$ the
functions $f\in\Dom(\boldsymbol{\D})$ satisfy certain
matching conditions. In order to describe them,
denote by $f_-$ the restriction $f\res e_y^-$ 
and by $f_j,\;j=1,\ldots,b(y)$
the restrictions $f\res e_y^j$.
The matching conditions at $y\neq o$ are
\begin{equation*}%\label{match}
f_-(y)=f_1(y)=\ldots=f_b(y);\ \
f_1'(y)+\ldots+f_b'(y)=f'_-(y).
\end{equation*}
The first condition comes
from the requirement $f\in \plainH1(\G)$
which includes continuity of $f$,
and the second appears as the
natural condition in the sense of
Calculus of Variations.
It is easy to check that the conditions listed are
also sufficient for $f\in\Dom (\boldsymbol{\D})$.

Let $V$ be a measurable, real-valued, bounded from below and \textsl{symmetric}
(that is, depending only on $|x|$) function on $\G$. 
Along with the Laplacian we
shall be interested also in the
Schr\"odinger operator with the potential $V(|x|)$:
\begin{equation}\label{lap:1}
\BA_Vf:=-\boldsymbol{\D} f+V(|x|)f.
\end{equation}
The operator $\BA_V$ is defined via its quadratic form 
\begin{equation}\label{lap:2}
\ba_V[f]:=\int_\G\bigl(|f'(x)|^2+V(|x|)|f(x)|^2\bigr)dx
\end{equation}
considered on the natural domain 
$\Dom(\ba_V)=\Quad(\CA_V)=\plainH{1,0}(\G)\cap \plainL2_V(\G)$.
On this domain the quadratic form \eqref{lap:2} is bounded from below
and closed in $\plainL2(\G)$, and the corresponding 
self-adjoint operator is taken as the realization of the operator
\eqref{lap:1}. We do not need the precise description of the domain
and the action of $\BA_V$. 

\section{The basic decomposition of $\plainL2(\G)$}
Consider two types of subtrees $T\subset\G$. Namely,
for any vertex $y$ and for any edge $e=\langle z,w\rangle$ we 
set 
\begin{equation*}%\label{tr:sub}
T_y=\{x\in\G:x\succeq y\},\qquad T_e=e\cup T_w.
\end{equation*}
Evidently $T_o=\G$. 

For any subtree $T=T_y$ or $T=T_e$ its
{\sl{branching function}} $g_T(t)$ 
is defined as
\begin{equation*}
g_T(t)=\#\{x\in T:|x|=t\}.
\end{equation*}
If $T=T_e$ and $\gen(e)=k\ge 0$, then $g_T(t)=g_k(t)$ where
\begin{equation*}
g_k(t)=\begin{cases} 0,\qquad  t<t_k,\\ 1,\qquad t_k\le t\le t_{k+1},\\
b_{k+1}\ldots b_n, \qquad t_n<t\le t_{n+1},\ n> k. \end{cases}
\end{equation*}
In particular, $g_0(t)=g_\G(t)=1$ for $0\le t\le t_1$ and 
\begin{equation*}
g_\G(t)=b_1\ldots b_n, \qquad t_n<t\le t_{n+1},\ n\ge 1.
\end{equation*}
So we see that 
\begin{equation}\label{tr:brf1}
g_k(t)=(b_1\ldots b_k)^{-1}g_\G(t),\qquad t>t_k,\ k>0.
\end{equation}
Note that 
\begin{equation}\label{tr:32}
\int_\G g_\G(t)dt=|\G|.
\end{equation}

Given a subtree $T\subset\G$, we say that a function
$f\in \plainL2(\G)$ belongs to the set
$\CF_T$ if and only if $f=0$ outside $T$
and
\begin{equation*}%\label{2}
f(x)=f(y)\qquad\text{if}\ x,y\in T\ \text{and}\ |x|=|y|.
\end{equation*}
Evidently $\CF_T$ is a closed subspace
of $\plainL2(\G)$. Any function $f\in \CF_{T_e}$, $\gen(e)=k\ge 0$ can be
naturally identified with the function $u:=J_e f$ on $(t_k,R)$, such that
$f(x)=u(|x|)$ for each $x\in T_e$ and $f(x)=0$ outside $T_e$. We have
\begin{gather}
\int_\G |f(x)|^2dx=\|u\|^2_{\plainL2((t_k,R),g_k)}:=
\int_{t_k}^R |u(t)|^2g_k(t)dt;\label{20}\\
f\in\CF_{T_e},\ u=J_e f.\notag
\end{gather}
This shows that the operator $J_e$ defines an isometry of the subspace
$\CF_{T_e}$ onto the weighted space $\plainL2\bigl((t_k,R);g_k\bigr)$. 
Along with \eqref{20}, we have
\begin{gather}
\int_\G |f'(x)|^2dx=a_k[u]:=
\int_{t_k}^R |u'(t)|^2g_k(t)dt;\label{21}\\
f\in\CF_{T_e}\cap\plainH{1,0}(\G),\ u=J_e f.\notag
\end{gather}

For the sake of brevity, 
below we use the notations
$\CF_y$ for $\CF_{T_y}$ and $\CF_y^j$ for $
\CF_{T_{e_y^j}},\ j=1,\ldots, b=b(y)$. 
It is clear that the subspaces
$\CF_y^1,\ldots,\CF_y^b$ are mutually
orthogonal and their orthogonal sum
\begin{equation*}%\label{lap:5}
\widetilde{\CF_y}=\CF_y^1\oplus\ldots\oplus\CF_y^b
\end{equation*}
contains $\CF_y$. Denote 
\begin{equation*}
\CF'_y = \widetilde{\CF_y}\ominus\CF_y 
\end{equation*}
\begin{theo}\label{lap:thm10} Let $\G$ be a regular tree.
\begin{itemize}
\item[(i)]
The subspaces $\CF'_y$, $y\in\CV(\G)$ are
mutually orthogonal and orthogonal to $\CF_\G$.
Moreover,
\begin{equation}\label{lap:10}
\plainL2(\G)=\CF_\G\oplus
\sum_{y\in\CV(\G)}\oplus\CF'_y.
\end{equation}
\item[(ii)]
Let $V(t)$ be a real, measurable and bounded below function on $\R_+$.
Then the decomposition \eqref{lap:10} reduces the Schr\"odinger
operator \eqref{lap:1}.
\end{itemize}
\end{theo}

According to Theorem 4.1, description of the spectrum
$\s(\BA_V)$
reduces to the similar
problem for the parts of $\BA_V$ in the
components of the decomposition
\eqref{lap:10}. These parts can be described in terms
of auxiliary differential operators  $\CA_{V,k}$,
$k=0,1,\ldots$ acting in the spaces $\plainL2\bigl((t_k,R), g_k\bigr)$.

For $f\in\Dom(\ba)\cap\CF_{T_e}$  the quadratic form
$\ba_V$, cf. \eqref{lap:2}, transforms as follows:
\begin{gather}
\ba_V[f]=a_{V,k}[u]:=\int_{t_k}^R\bigl(|u'(t)|^2+V(t)|u(t)|^2\bigr)g_k(t)dt,
\label{coa}\\ f\in\Dom(\ba)\cap\CF_{T_e},\ \gen(e)=k,\  u=J_e f.\notag
\end{gather}
For $V\equiv 0$ the quadratic form $a_{V,k}[u]$ turns into the
quadratic form $a_{k}[u]$ defined in \eqref{21}.
We define $\CA_{V,k}$ as the self-adjoint operator in 
$\plainL2\bigl((t_k,T), g_k\bigr)$, associated with the quadratic form
$a_{V,k}$.
We drop the subindex $V$ in these notation when dealing with the 
free Laplacian $-\boldsymbol{\D}=\BA_0$.

The following result was actually proved in \cite{C1} and \cite{NS2}; 
minor distinctions in the formulations are unessential.
\begin{theo}\label{thm110} Let $\G$ be a regular tree and $V$ be a bounded
from below, real-valued symmetric potential on $\G$. Then the part
of $\BA_V$ in the subspace $\CF_\G$ is unitarily equivalent
to the operator $\CA_{V,0}$ and the part of $\BA_V$ in each subspace
$\CF'_z$, $\gen(z)=k>0$, is unitary equivalent to the orthogonal sum
of $(b_k-1)$ copies of $\CA_{V,k}$.
\end{theo}
The next theorem is an immediate consequence of Theorems 4.1 and 4.2.
Below $\CA^{[r]}$ stands for the orthogonal sum of $r$
copies of a self-adjoint operator $\CA$. The symbol ``$\sim$'' means
unitary equivalence.

\begin{theo}\label{thm2} Under the assumptions of Theorem 4.2 the operator
$\BA_{V,\G}$ is unitary equivalent to the orthogonal sum of the operators
$\CA_{V,k}$, with growing multiplicities:
\begin{equation}\label{ortsum}
\BA_{V,\G}\sim \CA_{V,0}\oplus\sum_{k=1}^\infty\oplus 
{\CA_{V,k}}^{[b_1\ldots b_{k-1}(b_k-1)]}.
\end{equation}
\end{theo}

 \section{Main results} 
\subsection{The eigenvalue counting functions}\label{gennat}
For a self-adjoint, bounded from below operator $\CA$ with discrete 
spectrum,
we denote by $N(\l;\CA)$ the distribution function of its eigenvalues 
$\l_j(\CA)$ (counted according to their multiplicities),
\begin{equation*}
N(\l;\CA)=\#\{j:\l_j(\CA)<\l\},\qquad \l\in\R.
\end{equation*}

We start with the following simple but useful statement.
\begin{theo}\label{sp:11} Let $\G$ be a regular tree and
$V$ be a symmetric measurable real-valued function, bounded 
below. The spectrum 
$\s(\BA_V)$ is discrete if and only if the spectrum of the operator
$\CA_{V,0}$ is discrete. If this is the case, then
\begin{multline}\label{r:dec}
N(\l;\BA_V)=N(\l;\CA_{V,0})+\sum_{k=1}^\infty b_1\ldots b_{k-1}(b_k-1)
N(\l;\CA_{V,k}),\qquad\l\in\R.
\end{multline}
\end{theo}
{\bf Proof.} Consider the Rayleigh quotient 
$a_{V,k}[u]/\|u\|^2_{\plainL2((t_k,R),g_k)}$, cf. \eqref{20} and \eqref{coa}.
Due to the equality \eqref{tr:brf1}, this ratio does not change if we
replace in its numerator and denominator
the weight function $g_k(t)$ by $g_\G(t)$.
Now it follows from the variational 
principle that the spectrum of each operator $\CA_{V,k},\
k=1,2,\ldots$ is discrete
provided this is true for $k=0$. Moreover, we see that
\begin{equation*}
N(\l;\CA_{V,k_2})\le N(\l;\CA_{V,k_1}),\qquad k_1<k_2,\ \l>0. 
\end{equation*}
The discreteness of $\s(\CA_{V,k})$ for all $k$ implies the same property of 
$\s(\BA_V)$. The converse is evident. The equality \eqref{r:dec}
is an immediate consequence of the relation \eqref{ortsum}. 
\vskip0.2cm

The detailed study of the function $N(\l;\BA_V)$ is hampered by
the presence of the rapidly growing factors $b_1\ldots b_{k-1}(b_k-1)$.
These factors reflect geometry of the tree and 
do not depend on the potential $V$. For this reason, sometimes
we consider 
another counting function (introduced in \cite{SS}):
\begin{equation}\label{r:tilde}
\widetilde N(\l;\BA_V):=\sum_{k=0}^\infty N(\l;\CA_{V,k}).
\end{equation}

\subsection{The spectrum of the Laplacian}\label{splap} The
spectrum of the Laplacian on a regular tree depends on the behaviour of
the sequences $\{t_k\}$ and $\{b_k\}$ and can be quite different. We
present here several results in this direction. The proofs of those which
are new are given in the next section. In other cases we give the relevant
references. 
 
Our first result is quite elementary and its proof is standard. The result
applies to arbitrary metric graphs rather than to trees only. 
\begin{theo}\label{sp:1} Let $\G$ be a metric graph such that
$\sup_{e\in\CE(\G)}|e|=\infty$. Then the spectrum of the Laplacian
$-\boldsymbol{\D}$ on $\G$ coincides with $[0,\infty)$. 
\end{theo}
 
\bigskip
 
Other results concern the regular trees. We start with the trees of finite
radius, for which the information provided is rather complete. 

\begin{theo}\label{sp:4} 
Let $\G$ be a regular tree and $R(\G)<\infty$. 
\begin{itemize}
\item[(i)]
The spectrum of the Laplacian $-\boldsymbol{\D}$ on $\G$ is discrete. 
For each operator $\CA_{k}$  its eigenvalues
behave according to the 
Weyl law,
\begin{equation}\label{l:w}
\pi N(\l;\CA_{k})=\sqrt\l(R-t_k)+o(\sqrt\l),\qquad\l\to\infty.
\end{equation}
\item[(ii)]
If $|\G|<\infty$, 
then the Weyl asymptotic law
holds  for the operator $-\boldsymbol{\D}$:
\begin{equation}\label{l:ww}
\pi N(\l;-\boldsymbol{\D})=\sqrt\l|\G|+o(\sqrt\l),\qquad\l\to\infty.
\end{equation}
\item[(iii)]
If \begin{equation*}
\widetilde R(\G):=\sum\limits_{k=0}^\infty(R-t_k)<\infty,
\end{equation*}
then
\begin{equation}\label{l:wa}
\pi\widetilde N(\l;-\boldsymbol{\D})=\sqrt\l\widetilde R(\G)+o(\sqrt\l),
\qquad\l\to\infty.
\end{equation}
\end{itemize}
\end{theo}

Theorem 5.3 refines an earlier result of \cite{C2}. The case of general
(i.e. not necessarily regular) trees was analyzed in \cite{NS1}, Theorem 4.1. 
For the trees with $|\G|<\infty$, satisfying some additional assumptions,
the Weyl asymptotics \eqref{l:ww} follows from this theorem. However, the
result of \cite{NS1} does not cover the case of arbitrary regular trees 
of finite total length. 

The next statement can be derived from Theorem 5.3 by means of the elementary
variational arguments. We present it without proof.

\begin{cor}\label{sp:31} Let $\G$ be a regular tree, $R(\G)<\infty$, and
the potential $V(x)$ (not necessarily symmetric) be bounded. Then the
spectrum $\s(-\boldsymbol{\D}+V)$ is discrete.  If in addition
$|\G|<\infty$, then the asymptotic formula \eqref{l:ww} holds for its
eigenvalues. 

If in addition, the potential is symmetric, then the asymptotic formula
\eqref{l:w} holds for each operator $\CA_{V,k}$.  \end{cor}
 
\bigskip
 
If for a regular tree $\G$ one has $R(\G)<\infty$ but $|\G|=\infty$, then
the asymptotic behaviour of $N(\l;-\boldsymbol{\D})$ can be rather exotic.
The following example can serve as an illustration. 
 
Fix the numbers $q\in(0,1)$ and $b\in\N$.  Consider the
regular tree $\G=\G_{q,b}$ defined by the sequences $t_k=1-q^k$,
$k=0,1,\ldots$ and $b_k=b$, $k=1,\ldots$. Then $R(\G)=1$, so that the
spectrum of the Laplacian on $\G$ is always discrete. Further,
$g_0(t)=b^k$ for $t_k<t\le t_{k+1}$. The total length of $\G$ is
\begin{equation*} |\G|=1-q+\sum_{k=1}^\infty
b^k(q^k-q^{k+1})=(1-q)\sum_{k=0}^\infty(bq)^k.  
\end{equation*} Hence,
$|\G|=\frac{1-q}{1-bq}<\infty$ if $bq<1$ and $|\G|=\infty$ otherwise.  In
the first case, Theorem 5.3 (ii) shows that the Weyl law \eqref{l:ww} holds
for the eigenvalues of $-\boldsymbol{\D}$. Besides, 
$\widetilde R(\G_{q,b})=(1-q)^{-1}<\infty$, and by Theorem 5.3 (iii)
the asymptotic formula \eqref{l:wa} holds for any $q<1$ and any $b$.

Below we present the results for the function $N(\l,-\boldsymbol{\D})$, for 
$bq\ge1$. The case $bq>1$ was analyzed in \cite{NS1}, Example 8.2 (where
one should take $\a=0$). The result of \cite{NS1} for $bq=1$ was not
complete. 

\begin{theo}\label{sp:20} Let $\G=\G_{q,b}$. 
 \begin{itemize}
\item[(i)]
If $bq>1$, then there exists a bounded and bounded away from zero
periodic function $\psi$ with the period $\ln(q^{-2})$ such that
\begin{equation}\label{sp:45} \pi
N(\l;-\boldsymbol{\D})=\l^{\b/2}\bigl(\psi(\ln \l)+o(1)\bigr),\qquad
\l\to\infty \end{equation} where $\b=-\log_q b>1$. 
 \item[(ii)]
If $bq=1$, then \begin{equation}\label{sp:46} \pi
N(\l;-\boldsymbol{\D})=\frac{1-q}{2\ln b} \sqrt\l\bigl(\ln
\l+O(1)\bigr),\qquad \l\to\infty. 
\end{equation} 
\end{itemize}
\end{theo} 
The proof given
in the next section covers both cases. For $bq>1$, it reproduces the argument
from \cite{NS1}. 

\bigskip

The results for the trees with $R(\G)=\infty$ are much less exhaustive. 
We start with a criterion of positive definiteness of the Laplacian on
a regular tree, proven in \cite{NS2}. 
\begin{theo}\label{sp:35} Let $\G$ be
a regular tree and $R(\G)=\infty$. Then the Laplacian on $\G$ is positive
definite in  $\plainL2(\G)$ if and only if
\begin{equation}\label{sp:36}
\sup_{t>0}\biggl(\int_0^tg_\G(s)ds\cdot\int_t^\infty\frac{ds}{g_\G(s)}\biggr) 
<\infty. \end{equation} \end{theo}
 \vskip0.2cm
The condition \eqref{sp:36} is satisfied, in particular, for the
homogeneous trees $\G_b$. For them the spectrum can be described completely.
The next result is proven in \cite{SS}, Theorem 3.3. 
Introduce the number \begin{equation*}
\t=\arccos\frac{2}{b^{1/2}+b^{-1/2}}. \end{equation*}
 
\begin{theo}\label{sp:2} The spectrum of 
the operator $-\boldsymbol{\D}$ on the tree
$\G_b$ is of infinite multiplicity and consists of the bands
$\bigl[\bigl(\pi(l-1)+\t\bigr)^2,\bigl(\pi l-\t\bigr)^2\bigr]$, $l\in\N$
and the eigenvalues $\l_l=\bigl(\pi l\bigr)^2$. 
\end{theo} 
So, in this case
the spectrum has the {\sl band-gap structure} which is typical for periodic
problems. An analogue of Theorem 5.7 can be proved for regular trees for
which the sequences $t_{k+1}-t_k$ and $b_k$ are not necessarily
constant, as for $\G_b$, but periodic. 
 
\bigskip
 
Suppose now that for a regular tree all the branching numbers are equal,
$b_1=b_2=\ldots =b$, but the edge lengths $l_k=t_k-t_{k-1}$ are identically
distributed random variables. 
 
More precisely, let $[L_1,L_2]$ be a finite segment, $L_1>0$. Suppose that
$\mu$ is a Borelian probability measure on $[L_1,L_2]$. Denote by
$\mu^\infty$ the product of infinitely many copies of $\mu$; this is a
measure on the space of all sequences $\{l_k\}_{k\in\N}$ taking their
values in $[L_1,L_2]$. 
 
\begin{theo}\label{sp:3} Let the measure $\mu$ be absolute continuous. 
Suppose that for each $k=1,2,\ldots$ the lengths $l_k$ 
are independent random variables with distribution $\mu$.
Then
almost surely with respect to the measure $\mu^\infty$, the spectrum of
the operator $-\boldsymbol{\D}$ on $\G$ contains no absolute continuous
component. \end{theo}
 
The proof, which we do not present in this paper, was obtained in 
cooperation with G.Berkolaiko, K.Naimark, and U.Smilansky. Its starting
point is the equality \eqref{r:dec}. Then the spectrum of each operator
$\CA_k$ is analyzed with the help of F\"urstenberg's Theorem on the
product of random matrices.

Later the author had an opportunity to discuss this result with I.Goldsheid.
Here is the information provided by him.

1. The result (for the components $\CA_k$) was known to him and to S.Molchanov
before.

2. Moreover, the spectrum $\s(\CA_k)$ is almost surely pure point and
the eigenfunctions exponentially decay as $t\to\infty$ (the property which
is called Anderson localization).

\subsection{Operators $-\boldsymbol{\D}+V$ with growing potential}\label
{schr}

\begin{theo}\label{sp:32} Let $\G$ be a regular tree and $R(\G)=\infty$.
Denote by $\Psi$ the counting function for the sequence $\{t_k\}$,
\begin{equation*}
\Psi(\l)=\#\{k:t_k<\l\},\qquad\l>0.
\end{equation*}
Let $V(t)$ be a non-negative,
strictly monotonically increasing, unbounded
continuous function on $\R_+$. Let $Q$ stand for its inverse. Suppose
that the functions $Q$ and $\Psi\circ Q$ satisfy the $\D_2$-condition
\begin{equation}\label{del} 
Q(2\l)\le CQ(\l),\qquad \l\ge\l_0; 
\end{equation} 
\begin{equation}\label{del1} \Psi(Q(2\l))\le C\Psi(Q(\l)),\qquad
\l\ge\l_0 \end{equation} and that 
\begin{equation}\label{ma}
\Psi(t)=o\bigl(t\sqrt{V(t)}\bigr),\qquad t\to\infty. 
\end{equation} Then for the counting function $\widetilde N(\l;\BA_V)$,
cf. \eqref{r:tilde}, the asymptotic formula is valid:
\begin{equation}\label{sch} \pi \widetilde
N(\l;\BA_V)=(1+o(1))\sum_{k=0}^\infty\int_{t_k}^\infty
(\l-V(t))_+^{1/2}dt,\qquad\l\to\infty. 
\end{equation}

\end{theo}

 The asymptotic formula \eqref{sch} looks
quite natural. The condition \eqref{del} is standard for this class of 
problems. Two other conditions are rather restrictive and we do not know
whether they are sharp.
Note that the condition \eqref{del1} is automatically satisfied if 
we suppose that the function $\Psi$ itself satisfies the $\D_2$-condition.
Note also that for $t_k=k^r,\ r>0$ and $V(t)=t^\g,\ \g>0$ the assumption
\eqref{ma} reduces to $r^{-1}<1+\g/2$.

\section{Proofs}
\subsection{ Proof of Theorem 5.2}
It is enough to show that for any $r>0$ the point $\l=r^2$ belongs to 
$\s(-\boldsymbol{\D})$.
For this purpose we fix a non-negative function $
\varf\in C^\infty_0(-1,1)$ such that $\varf(t)=1$ on $(-1/2,1/2)$.
Further, choose an edge $e\in\CE(\G)$. In an appropriate coordinate
system, $e$ can be identified with the interval $(-l,l)$
where $l=|e|/2$. The function $f$ on $\G$,
\begin{equation*}
f(t)=\varf(t/l)\sin rt\ {\text {on}}\ e,\qquad  f(t)=0\ {\text {otherwise}},
\end{equation*}
belongs to $\Dom(\boldsymbol{\D})$.
An elementary calculation shows that
\begin{equation*}
\|\boldsymbol{\D} f+r^2f\|\le\vare(l)\|f\|,\qquad \vare(l)\to 0\ 
{\text {as}} \ l\to\infty.
\end{equation*} 
Choosing a sequence of edges $e$ such that $|e|\to\infty$, we obtain
a Weyl sequence for the operator $-\boldsymbol{\D}$ and the point $\l=r^2$.
This implies that $\l\in\s((\boldsymbol{\D})$.

\bigskip

\subsection{ Auxiliary material} 
We shall use the variational techniques, in the spirit of the book \cite{BS}.
We present the material we need in the form, convenient for the applications
to the operators $\CA_k$.

Let $w(t)$ be a monotonically growing function on a finite
interval $[a,b)$. In our applications we shall take $[a,b)=[t_k,R)$ and
$w(t)=g_k(t)$, which explains the nature of 
our assumptions about the function $w$.
We suppose that $w(t)\ge1$ and that the points $t_k$ of 
discontinuity of $w$ may accumulate at the point $b$ only.
Consider the Hilbert space $\CH^{1,\bullet}\bigl((a,b),w\bigr)$ whose 
elements are the functions $u$ on $[a,b)$, such that 
$u\in\plainH1(a,b-\e)$ for any $\e>0$, $u(a)=0$, and
\begin{equation*}
\|u\|^2_{\CH^{1,\bullet}((a,b)),w)}:=\int_a^b|u'(t)|^2w(t)dt<\infty,
\end{equation*} 
cf. \eqref{21}. 
We write $\CH^{1,\bullet}(a,b)$ instead of 
$\CH^{1,\bullet}\bigl((a,b),1\bigr)$.
The weighted Sobolev space with the weight $w$ is defined as 
\begin{equation*}
\plainH{1,\bullet}\bigl((a,b),w\bigr) =
\CH^{1,\bullet}\bigl((a,b),w\bigr)\cap\plainL2\bigl((a,b),w\bigr).
\end{equation*}

Let us change the variables, taking
\begin{equation}\label{pr:ch}
s=s(t)=\int_a^t\frac{d\tau}{w(\tau)}.
\end{equation} 
The variable $s$ runs over the interval $[0,L)$ where
\begin{equation*}
L=\int_a^b\frac{d\tau}{w(\tau)}.
\end{equation*}
Since $w(t)\ge1$, we have $L\le b-a<\infty$. 
Below $t(s)$ stands for the function on $[0,L)$, inverse to $s(t)$. The
derivative $t'(s)=w(t(s))$ exists everywhere, 
except for the points $s_k=s(t_k)$.

Let $y(s)=u(t(s))$, then
\begin{equation}\label{pr:0}
\|u\|^2_{\CH^{1,\bullet}((a,b),w)}=\int_0^L|y'(s)|^2ds
\end{equation} 
and
\begin{equation*}%\label{pr:00}
\|u\|^2_{\plainL2((a,b),w)}=\int_0^LW(s)|y(s)|^2ds,
\qquad W(s)=w^2(t(s)).
\end{equation*} 
%It follows from \eqref{pr:0} that for any 
%$u\in\CH^{1,\bullet}\bigl((a,b),w\bigr)$ its image $y$ has a finite limit
%at $s=L$. Therefore, the same is true for the function $u$ at the point
%$t=b$.

The function $W(s)$ is monotone, and
\begin{equation}\label{pr:1}
\int_0^L W(s)ds=\int_0^L w(t(s))t'(s)ds=\int_a^b w(t)dt;
\end{equation}
\begin{equation}\label{pr:2}
\int_0^L \sqrt{W(s)}ds=\int_0^L t'(s)ds=b-a.
\end{equation}

\bigskip

In the course of the proofs of Theorems 5.3 and 5.5
we make use of the following
result. Its most important part (i) 
was obtained in \cite{BB}, see Theorem 3.1 and,
especially, Remark 3.1 there.
See also an exposition in \cite{BLS}, Corollary 6.3. The part (ii) is new
and we present it with proof.

\begin{theo}\label{sp:40} 
\begin{itemize}
\item[(i)] Let $L\le\infty$ and let $W\in\plainL{1/2}(0,L)$ 
be a monotone, non-negative function. Then the inequality holds
\begin{equation}\label{pr:5}
\int_0^L W(s)|y(s)|^2ds\le C(W)
\int_0^L |y'(s)|^2ds,\qquad y\in\CH^{1,\bullet}(0,L),
\end{equation}
and therefore the quadratic form in the left-hand side
generates in $\CH^{1,\bullet}(0,L)$ a bounded self-adjoint operator,
say $T_W$. Moreover, the operator $T_W$
is compact and for its eigenvalues $\mu_j(T_W)$ the following estimate holds,
with a constant factor which  does not depend on $L$ and on $W$:
\begin{equation}\label{l:est}
\#\{j:\mu_j(T_W)>\l^{-1}\}\le C\sqrt\l \int_0^L\sqrt{W(s)}ds,\qquad \l>0.
\end{equation}
Also, the 
asymptotic formula is valid:
\begin{multline}\label{l:we}
\#\{j:\mu_j(T_W)>\l^{-1}\} =\frac{\sqrt\l}{\pi}\int_0^L\sqrt{W(s)}ds+
o(\sqrt\l),
\qquad\l\to\infty.
\end{multline}
\item[(ii)] Suppose in addition that the function $W$ satisfies the estimate
\begin{equation}\label{l:est2}
W(s)\le C(L-s)^{-r}
\end{equation}
with some $r\in(0,2)$. Then the following, uniform in $\l$ remainder estimate
in the asymptotic formula \eqref{l:we} is satisfied:
\begin{multline}\label{l:weres}
\biggl|\#\{j:\mu_j(T_W)>\l^{-1}\}-\frac{\sqrt\l}{\pi}\int_0^L\sqrt{W(s)}ds
\biggr|\\ \le C(L)\bigl(\l^{1/(4-r)}+1\bigr),\qquad\l>0.
\end{multline}
\end{itemize}
\end{theo}
{\bf Proof of (ii).} For definiteness,
we assume the function $W$ to be increasing. 

Suppose at first that $W$ is bounded. Then we use the
standard variational reasoning: divide $[0,L]$ into $n$ equal parts,
on each part $(s_k,s_{k+1})$ replace $W(s)$ by its $\inf$ and $\sup$ and
solve the resulting eigenvalue problem under the Dirichlet or the Neumann
boundary conditions. We obtain, denoting $h=L/n$:
\begin{equation*}
\sum_{k=0}^{n-1}\biggl[\frac{h}{\pi}\sqrt{\l W(s_k+)}\biggr]
\le \#\{j:\mu_j(T_W)>\l^{-1}\}\le
n +\sum_{k=1}^{n} \biggl[\frac{h}{\pi}\sqrt{\l W(s_k-)}\biggr]. 
 \end{equation*}
Roughening this inequality, we obtain:
\begin{multline*}%\label{2:4}
-n+\frac{h\sqrt\l}{\pi}\sum_{k=0}^{n-1}\sqrt{ W(s_k+)}
\le  \#\{j:\mu_j(T_W)>\l^{-1}\}\\ \le
n +\frac{h\sqrt\l}{\pi}\sum_{k=1}^{n} \sqrt{ W(s_k-)}. 
 \end{multline*}
We also have, due to the monotonicity of $W$:
\begin{equation*}
 h\sum_{k=0}^{n-1}\sqrt{W(s_k+)}\le\int_0^L \sqrt{W(s)}ds\le 
h\sum_{k=1}^{n} \sqrt{W(s_k-)}.
\end{equation*}
This yields
\begin{multline}\label{2:5}
\biggl| \#\{j:\mu_j(T_W)>\l^{-1}\}-\frac{\sqrt{\l}}
{\pi}\int_0^L \sqrt{W(t)}dt\biggr|
\le
C\bigl(\sqrt{\l}\frac{L\sqrt{W(L-)}}{n}+n\bigr).
\end{multline}

Suppose now that $W(s)$ is unbounded. Then we choose a point $S<L$,
insert the condition $y(S)=0$,
apply the inequality \eqref{2:5} on $(0,S)$ and use the estimate \eqref{l:est}
on $(S,L)$. We obtain
\begin{multline*}
\biggl|\pi\#\{j:\mu_j(T_W)>\l^{-1}\} -\sqrt{\l}\int_0^L
\sqrt{W(s)}ds\biggr|\\ \le
C\biggl(\sqrt{\l}\biggl(\frac{L\sqrt{W(S)}}{n}+\int_S^L \sqrt{W(s)}ds\biggr)+
n+1\biggr).
\end{multline*}
Now we use the inequality \eqref{l:est2} and 
then minimize the right-hand side over $S\in(0,L)$. This gives
\begin{multline*}
\biggl|\pi\#\{j:\mu_j(T_W)>\l^{-1}\} -\sqrt{\l}\int_0^L \sqrt{W(s)}ds\biggr|
 \le
C\biggl(\sqrt{\l}(L/n)^{1-(r/2)}+n+1\biggr).
\end{multline*}

We arrive at \eqref{l:weres}, taking here 
$n=\bigl[\l^{1/(4-r)}\bigr]+1$.

\subsection{ Proof of Theorem 5.3}
(i) Fix $k=0,1,\ldots$ and apply the construction in the beginning
of Subsection 6.2 to the interval
$[a,b)=[t_k,R)$ and the weight function $w(t)=g_k(t)$. Let $s_k(t)$
stands for the corresponding function \eqref{pr:ch} and $t_k(s)$
stands for its inverse.
The assumptions
of Theorem 6.1 are satisfied for the function $W_k(s)=g_k^2(t_k(s))$.
According to \eqref{pr:2}, the relations \eqref{l:est} 
and \eqref{l:we}
turn into
\begin{equation}\label{l:est1}
\#\{j:\mu_j(T_{W_k})>\l^{-1}\}\le C\sqrt\l (R-t_k),\qquad \l>0,
\end{equation}
where the constant $C$ does not depend on $k$, and
\begin{equation}\label{l:we1}
\#\{j:\mu_j(T_{W_k})>\l^{-1}\}  =\frac{\sqrt\l}{\pi}(R-t_k)+o(\sqrt\l),
\qquad\l\to\infty.
\end{equation}

Substituting $u(t)=y(s_k(t))$
in \eqref{pr:5} (with $L=L_k=\int_{t_k}^R (g_k(\tau))^{-1}d\tau$
and $W=W_k$), we come to the inequality
\begin{equation*}
\int_{t_k}^R|u(t)|^2g_k(t)dt\le C a_k[u],
\qquad u\in\CH^{1,\bullet}\bigl((t_k,R),g_k\bigr)
\end{equation*}
where $a_k[u]$ is the quadratic form defined in \eqref{21}.
This shows that the operator $\CA_k$ has bounded inverse and
that the spectrum of $\CA_k^{-1}$ coincides with the one of the
operator $T_{W_k}$. The
inequality \eqref{l:est1} turns into

\begin{equation}\label{l:ek}
N(\l;\CA_k)\le C\sqrt\l (R-t_k),\qquad \l>0, k=0,1,\ldots
\end{equation} and
the asymptotic formula \eqref{l:we1} turns into the formula \eqref{l:w}.
\vskip0.2cm
(ii) By \eqref{r:dec}, we have
\begin{multline*}
\l^{-1/2} N(\l;-\boldsymbol{\D})=\biggl(\l^{-1/2} N(\l;\CA_{0})\biggr)\\
+\sum_{k=1}^\infty b_1\ldots b_{k-1}(b_k-1)
\biggl(\l^{-1/2} N(\l;\CA_{k})\biggr).
\end{multline*} 
As $\l\to\infty$, each term in big parentheses tends to $\pi^{-1}(R-t_k)$. 
Besides,
by \eqref{l:ek} the series is dominated by 
\begin{equation*}
C\bigl(R+\sum_{k=1}^\infty b_1\ldots b_{k-1}(b_k-1)(R-t_k)\bigr)=C|\G|.
\end{equation*}
Now \eqref{l:ww} follows from the Lebesgue Theorem on the 
dominated convergence.
\vskip0.2cm
(iii) The proof of \eqref{l:wa} is the same and we skip it.
\vskip0.2cm

\noindent{\bf Remark.}
It follows from \eqref{pr:0} that for any 
$u\in\CH^{1,\bullet}\bigl((t_k,R),g_k\bigr)$ its image $y$ has a finite limit
at $s=L_k$. Therefore, the same is true for the function $u$ at the point
$t=R$. The equalities \eqref{pr:1} and \eqref{tr:32} imply that necessarily 
$u(R-)=0$ 
for any function
$u\in\Dom(a_k)$, provided $|\G|=\infty$.
If $|\G|<\infty$, then various boundary conditions at $t=R$ for functions
are possible. This is
consistent with the result of \cite{C2}, Theorem 5.2 where the boundary
value problems for the differential equations on regular trees were studied
from a different point of view.

\subsection{Proof of Theorem 5.5} 
For the tree $\G_{q,b}$ we have for $k=0,1,\ldots$, 
cf. \eqref{20} and \eqref{21}:
\begin{equation*}
\|u\|^2_{\plainL2((t_k,R);g_k)}=\int_{1-q^k}^1 |u(t)|^2g_k(t)dt
\end{equation*}
and 
\begin{equation*}
a_k[u]=\int_{1-q^k}^1 |u'(t)|^2g_k(t)dt.
\end{equation*}
Substituting $t=1-q^k(1-s)$, $u(t)=v(s)$, and taking \eqref{tr:brf1}
into account, we obtain:
\begin{equation*}
\|u\|^2_{\plainL2((t_k,R),g_k)}=q^{-k}b^k\|v\|^2_{\plainL2((0,R),g_\G)};
\qquad a_k[u]=q^{k}b^ka_0[v].
\end{equation*}
This implies that for any $k$ the operator $\CA_k$ is unitarily equivalent
to $q^{-2k}\CA_0$ and therefore,

\begin{equation}\label{2:9}
N(\l;\CA_k)=N(\l q^{2k};\CA_0),\qquad \l>0,\ k=1,2,\ldots
\end{equation}
This property of self-similarity is the key observation which allows us
to handle the problem.

It follows from Theorem 4.3 and the formula \eqref{2:9} that
\begin{equation}\label{2:10}
N(\l;-\boldsymbol{\D})=N(\l;\CA_0)+(1-b^{-1})
\sum_{k=1}^\infty b^k N(\l q^{2k};\CA_0).
\end{equation}

Now we are in a position to complete the proof of the statement (i).
The function $N(\l;\CA_0)$ satisfies the inequality
\begin{equation}\label{n}
N(\l;\CA_0)\le C\sqrt\l,\qquad \l>0,
\end{equation}
 cf. \eqref{l:ek}. 
Denote $\mu=\ln\l$, $\eta=-2\ln q$, 
$\Phi(\mu)=\l^{-\b/2}N(\l;-\boldsymbol{\D})$ and 
$\varf(\mu)=\l^{-\b/2}N(\l;\CA_0)$, then the equality \eqref{2:10} turns into
\begin{equation*}
\Phi(\mu)=\varf(\mu)+(1-b^{-1})\sum_{k=1}^\infty\varf(\mu-k\eta ).
\end{equation*}
This yields
\begin{equation}\label{renew}
\Phi(\mu)-\Phi(\mu-\eta)=\varf(\mu)-b^{-1}\varf(\mu-\eta).
\end{equation}
This is 
a particular case of the {\sl Renewal Equation}, well known in probability. 
The function
in the right-hand side of \eqref{renew} is zero at $-\infty$ (since 
$N(\l;\CA_0)=0$ for small $\l>0$) and exponentially decays at $+\infty$
(since $N(\l;\CA_0)$ satisfies \eqref{n} and $\b>1$). Therefore, the
Renewal Theorem applies, see e.g. \cite{F}, Chapter XI.1, or a modern 
exposition in \cite{LV}. The equation \eqref{renew} involves the single
shift (by $\eta$), hence this is the so-called lattice case. According
to the Renewal Theorem, there exists an $\eta$-periodic function $\psi(\mu)$
which is bounded and bounded away from zero, such that
\begin{equation*}
\Phi(\mu)=\psi(\mu)+o(1),\qquad\mu\to\infty.
\end{equation*}
This immediately leads to \eqref{sp:45}.

(ii) As in the proof of Theorem 5.3, we use the scheme presented in the 
beginning of Subsect. 6.2. For the tree $\G=\G_{b^{-1},b}$  and $k=0$ we have 
$[a,b)=[0,1)$ and
\begin{equation*}
w(t)=g_\G(t)=b^j\ {\text for}\ 1-b^{-j}< t\le 1-b^{-j-1},\ j=0,1,\ldots.
\end{equation*}
It follows that
\begin{equation}\label{l:l}
w(t)\le C(1-t)^{-1}.
\end{equation}
Besides,
\begin{equation*}
L=\int_0^1\frac{dt}{w(t)}=\sum_{j=0}^\infty\frac{b^{-j}-b^{-j-1}}{b^j}=
\frac{b}{b+1}.
\end{equation*}
For the function $s=s(t)$ defined by \eqref{pr:ch}, we have
\begin{equation*}
s(t_k)=\sum_{j=0}^{k-1}\frac{b^{-j}-b^{-j-1}}{b^j}=L-\frac{b^{-2k}}{b+1},
\end{equation*}
or
\begin{equation*}
L-s(t_k)=\frac{(1-t_k)^2}{b+1}.
\end{equation*}
Since $s(t)$ is monotone, this implies
\begin{equation*}
L-s(t)\ge c(1-t)^2,\ c>0,\qquad t\in(0,1)
\end{equation*}
and therefore, $1-t(s)\le c^{-1/2}(L-s)^{1/2}$.
It follows from here and from \eqref{l:l} that 
\begin{equation*}
W(s)=w^2(t(s))\le C(L-s)^{-1},
\end{equation*}
so that the inequality \eqref{l:est2} is satisfied with $r=1$. 
Correspondingly, the estimate \eqref{l:weres} takes the form
\begin{equation}\label{l:res}
\bigl|\pi N(\l;\CA_0)-\l^{1/2}\bigr|\le C(\l^{1/3}+1),\qquad\l>0.
\end{equation}
\bigskip
Let us return to the equality \eqref{2:10} (where now $q=b^{-1}$).
According to the estimate \eqref{l:ek}, with $R=1$ and $k=0$, we find that
\begin{equation*}%\label{l:zer}
N(\l;\CA_0)=0\qquad{\text {if}}\ C^2\l<1.
\end{equation*}
Therefore, summation in \eqref{2:10} is actually taken over such $k$
that $b^{2k}\le C^2\l$. Using for all such $k$ the estimate 
\eqref{l:res}, we obtain:
\begin{multline*}
\biggl|\pi N(\l;\CA_0)-\biggl(1+(1-q)
\#\{k>1:\ b^{2k}\le C^2\l\}\biggr)\l^{1/2}
\biggr|\\ \le C
\sum_{k:b^{2k}\le C^2\l}\bigl(\l^{1/3}b^{k/3}+1\bigl).
\end{multline*}
The sum in the right-hand side is of order $O(\l^{1/2})$, and the factor
in front of $\l^{1/2}$ in the left-hand side 
differs from $\frac{(1-q)\ln \l}{2\ln b}$ by $O(1)$. This
completes the proof of
\eqref{sp:46}.

\subsection{Proof of Theorem 5.9} We split the proof into several steps.

{\bf 1.} Denote 
\begin{equation*}%\label{pr:j}
J(\l;V)=\int_0^\infty(\l-V(t))_+^{1/2}dt=\int_0^{Q(\l)} (\l-V(t))^{1/2}dt.
\end{equation*} 
Under the assumption \eqref{del} one has 
\begin{equation*}%\label{d} 
J(\l;q)\asymp Q(\l)\sqrt\l,\qquad\l\to\infty 
\end{equation*}
where the symbol $\asymp$ means a two-sided estimate.
Indeed, evidently
$J(\l;V)\le Q(\l)\sqrt\l$, and for $\l\ge 2\l_0$ one has 
\begin{equation*}
J(\l;V)=\frac{1}{2}\int_{V(0)}^\l \frac{Q(s)ds}{(\l-s)^{1/2}}
\ge\frac{Q(\l/2)}{2}\int_{\l/2}^\l \frac{ds}{(\l-s)^{1/2}}\ge
cQ(\l)\sqrt\l. 
\end{equation*}
Later we shall need also the inequality
\begin{equation}\label{pr:r}
\int_r^\infty(\l-V(t))_+^{1/2}dt\ge cQ(\l)\sqrt\l,\qquad r\le Q(\l/2),
\ \l\ge 2\l_0.
\end{equation}
Its proof, and also the value of $c$, are the same as in the preceding
inequality.
\vskip0.2cm
{\bf 2.} Consider the Schr\"odinger operator $K_Vy=-y''+Vy,\ u(y)=0$ in
$\plainL2(\R_+)$. Fix $\hat\l>V(0)$ and compare the values of $N(\hat\l;K_V)$
and $\#\{j:\mu_j(T_W)>1\}$ for the operator $T_W$ introduced in Theorem 6.1,
with $L=Q(\hat\l)$ and $W(s)=\hat\l-V(s)$. It follows from the decoupling
principle and from the Birman -- Schwinger principle that these two numbers
differ no more than by $2$. For estimating the number $\#\{j:\mu_j(T_W)>1\}$,
we use the inequality \eqref{2:5} (with $\l=1$). The only difference
with \eqref{2:5} is that now the function $W(s)$ is decreasing, and for
this reason the term $W(L-)$ in the right-hand side must be replaced by
$W(0)$. As a result, we obtain (replacing in the result $\hat\l$ by $\l$):

\begin{multline*}
\biggl|\pi N(\l;K_V)-J(\l;V)\biggr|\le 
C\biggl(\frac{Q(\l)\sqrt{\l-V(0)}}{n}+n+1\biggr)\\
\le C\biggl(\frac{Q(\l)\sqrt{\l}}{n}+n+1\biggr).
\end{multline*}
Let us stress that the factor $C$ does not depend on
$\l$ and $n$. Minimizing the right-hand side of the last inequality over $n$,
we come to the estimate
\begin{equation}\label{aa}
\biggl|\pi N(\l;K_V)-J(\l;V)\biggr|\le 
C\biggl(\bigl(Q(\l)\sqrt{\l}\bigr)^{1/2}+1\biggr)
\end{equation}
where $C$ is an absolute constant.
\vskip0.2cm
{\bf 3.} Consider the quadratic forms  $a_{V,k}$ defined in \eqref{coa}.
The corresponding operators $\CA_{V,k}$ act in the weighted spaces
$\plainL2\bigl((t_k,\infty),g_k\bigr)$ (recall that in our case $R=\infty)$.
For us it is more convenient to deal with the operators acting in the 
``usual'' $\plainL2$. For this purpose we make a substitution which we 
describe for $k=0$. The changes needed in the  case $k>0$, are evident.

Denote $y(t)=\sqrt{g_\G(t)}u(t)$, then also 
$y'(t)=\sqrt{g_\G(t)}u'(t),\ t\neq t_k,\ k=1,2,\ldots$, since $g_\G$ is
a step function. Evidently, 
\begin{equation*}
 \int_{\R_+}|u(t)|^2 g_\G(t)dt=\int_{\R_+}|y(t)|^2dt
\end{equation*}
and
\begin{equation}\label{pr:m}
a_{V,0}[u]=\widetilde{a}_{V,0}[y]:=\sum_{k=0}^\infty\int_{t_k}^{t_{k+1}}
\bigl(|y'(t)|^2+V(t)|y(t)|^2\bigr)dt.
\end{equation}
The domain $\Dom(a_{V,0})$ consists of all functions $y(t)$, 
such that $y\restriction(t_k,t_{k+1})\in\plainH1(t_k,t_{k+1})$ for each $k$,
the sum in the last side of \eqref{pr:m} is finite, $y(0)=0$ and
the matching conditions at the points $t_k$ are fulfilled:
\begin{equation*}
y(t_k+)=\sqrt{b_k}y(t_k-),\qquad k=1,2,\ldots
\end{equation*}

It is not difficult to derive from \eqref{aa} a similar estimate for the
operator $\CA_{V,0}$. Indeed, again the problem reduces to the interval
$(0,Q(\l))$. The linear spaces $\Dom(a_{V,0})$ and  $\Quad(K_V)$,
restricted to the set of functions supported by this interval, 
differ by a subspace of dimension $2\Psi(Q(\l))$. Therefore, \eqref{aa}
implies
\begin{equation*}
\biggl|\pi N(\l;\CA_{V,0})-J(\l;V)\biggr|\le 
C\biggl(\bigl(Q(\l)\sqrt{\l}\bigr)^{1/2}+\Psi(Q(\l))\biggr).
\end{equation*}
The term $1$ appearing in \eqref{aa} can be dropped, since $\Psi(Q(\l))\ge1$
for any $\l>V(0)$.

Quite similarly,
\begin{equation*}%\label{ff}
\biggl|\pi N(\l;\CA_{V,k})-\int_{t_k}^\infty(\l-V(t))_+^{1/2}\biggr|\le 
C\biggl(\bigl(Q(\l)\sqrt{\l}\bigr)^{1/2}+\Psi(Q(\l))\biggr).
\end{equation*}
\vskip0.2cm
{\bf 4.} Now we are in a position to complete the proof. We have
\begin{equation*}
\widetilde N(\l;\BA_V)=\sum_{k=0}^\infty N(\l;\CA_{V,k})
=\sum_{k:0\le t_k<Q(\l)}N(\l;\CA_{V,k}).
\end{equation*}
Therefore, 
\begin{multline}\label{pp}
\biggl|\widetilde N(\l;\BA_V)-\sum_{k=0}^\infty
\int_{t_k}^\infty(\l-V(t))_+^{1/2}dt\biggr|\\ \le C\Psi(Q(\l))
\biggl(\bigl(Q(\l)\sqrt{\l}\bigr)^{1/2}+\Psi(Q(\l))\biggr).
\end{multline}

Our next task is to estimate from below the sum appearing in the left-hand side
of \eqref{pp}. We derive from \eqref{pr:r}: 
\begin{multline*}
\sum_{k=0}^\infty\int_{t_k}^\infty(\l-V(t))_+^{1/2}dt\ge
\sum_{k:t_k\le Q(\l/2)}\int_{t_k}^\infty(\l-V(t))_+^{1/2}dt\\ \ge
c\Psi(Q(\l/2))Q(\l)\sqrt\l\ge c'\Psi(Q(\l))Q(\l)\sqrt\l.
\end{multline*}
The latter inequality is implied by the $\D_2$-condition \eqref{del1}.

It follows from the assumption \eqref{ma} that
\begin{equation*}
\Psi(Q(\l))\biggl(\bigl(Q(\l)\sqrt{\l}\bigr)^{1/2}+\Psi(Q(\l))\biggr)
=0\biggl(\Psi(Q(\l))Q(\l)\sqrt\l\biggr).
\end{equation*}
The desired asymptotic formula \eqref{sch} immediately follows.

\section{Regular trees without boundary}
In conclusion, let us discuss the case when the tree $\G$ has no boundary.

Let $\G$ be a general metric tree. Choose a vertex $o\in\CE(\G)$ and
suppose that there are $d$ edges of $\G$ adjacent to $o$. Then $\G$ can
be split into $d$ rooted subtrees $\G^1,\ldots\G^d$ having the 
common root $o$. We say that the tree $\G$ is regular if 
and only if all the subtrees
$\G^j$ are regular in the sense of Definition 2.1 and the
corresponding sequences $\{t_k\}$ and $\{b_k\}$ are the same for all 
$j=1,\ldots,d$. Note that this definition is not invariant with respect
to the choice of the vertex $o$.

Suppose now that all the subtrees $\{\G^j\}$ are homogeneous,
$\G^1=\ldots=\G^d=\G_b$ and $d=b+1$. Then we say that the tree $\G$ is
homogeneous. Unlike the case of arbitrary regular trees, this definition
is invariant with respect
to the choice of $o$.

The definitions of the Laplacian and of the Schr\"odinger operator 
extend to the trees without boundary in a natural way. The only difference
is that now we have no boundary condition at $o$. Instead, the functions from
the quadratic domain of the operator are required to be continuous at $o$.

Replacement of this continuity condition by the Dirichlet boundary condition
$u(o)=0$ means the passage to a subspace of codimension $1$ of $\Quad(\BA_V)$.
Therefore, the character of the spectrum is not affected, and moreover,
the eigenvalue distribution function $N(\l;\BA_V)$ can change no more than by
one. The  new operator splits into
the orthogonal sum of $d$ copies of the operator studied in the main part
of this paper. This allows one to immediately reformulate all the results
for this new situation. It is unnecessary to 
present their precise formulations. 

Note that the papers \cite{C1} and \cite{C2} deal with the operators
on trees without boundary.

\section{Acknowledgements}

The work on this paper was supported by the Minerva center for non-linear 
physics
and by the Israel Science Foundation.

I use this opportunity to express my deep thanks to G. Berkolaiko,
K. Naimark and U. Smilansky for their participation in the proof of
Theorem 10, and to I. Goldsheid for the useful discussion of this result.
My special gratitude goes
to U. Smilansky for his enthusiastic support of the study 
of differential operators on metric graphs.

%%%%%%%%%%%%%%%%%%%%%%%%%%%%%%%%%%%%%%%%%%%%%%%%%%%%%%%%%%%%%%%%%%%%%%%
%	Please insert your list of references below
%%%%%%%%%%%%%%%%%%%%%%%%%%%%%%%%%%%%%%%%%%%%%%%%%%%%%%%%%%%%%%%%%%%%%%%

\vspace*{\fill}
\noindent
Michael Solomyak \\

\noindent
Department of Mathematics\\

\noindent
The Weizmann Institute of Science\\

\noindent
Rehovot 76100\\

\noindent
Israel

\noindent
solom@wisdom.weizmann.ac.il
\end{document}